\documentclass{article}
\usepackage{amsmath,amssymb}
\title{Approximating the parallel transport of an induced connection}
\author{Derek Harland}
\date{4th April 2022}

\begin{document}

%Second order method:
%\begin{equation}
%\Omega^{(2)}(h,0)=\frac{1}{2}v^\dagger(h)v(0)+\frac{1}{2}(v^\dagger(0)v(h))^{-1}
%\end{equation}
%Third order method:
%\begin{equation}
%\Omega^{(3)}(h,0)=\frac{1}{3}\big[ 4\Omega^{(2)}\left(h,\tfrac{h}{2}\right)\Omega^{(2)}\left(\tfrac{h}{2},0\right) - \Omega^{(2)}(h,0) \big]
%\end{equation}
%Fourth order method:
%\begin{multline}
%\Omega^{(4)}(h,0)=\frac{1}{44}\big[
%90\Omega^{(2)}(h,\tfrac{2h}{3})\Omega^{(2)}(\tfrac{2h}{3},\tfrac{h}{3})\Omega^{(2)}(\tfrac{h}{3},0)\\
%-27\Omega^{(2)}(h,\tfrac{2h}{3})\Omega^{(2)}(\tfrac{2h}{3},0)
%-27\Omega^{(2)}(h,\tfrac{h}{3})\Omega^{(2)}(\tfrac{h}{3},0)
%+8\Omega^{(2)}(h,0)\big].
%\end{multline}
\maketitle
\abstract{Efficient numerical methods to approximate the parallel transport operators of the induced connection on a sub-bundle of a vector bundle are presented.  These methods are simpler than naive applications of a Runge--Kutta algorithm, and have accuracy up to order 4.  They have the desirable property of being insensitive to choices of trivialisation of the sub-bundle.  The methods were developed in order to solve a problem of computing skyrmions using the Atiyah--Manton--Sutcliffe and Atiyah--Drinfeld--Hitchin--Manin constructions, but are applicable to a broader range of problems in computational geometry.}

\section{Introduction}
\label{sec:1}
Given a hermitian vector bundle equipped with a unitary connection $\nabla$, any sub-bundle $E$ comes equipped with a natural connection $\nabla^E$.  The covariant derivative $\nabla^E s$ of any section $s$ of $E$ is defined to be the orthogonal projection of $\nabla s$ onto $E$.  Induced connections feature prominently in submanifold geometry, where the tangent and normal bundles of a submanifold inherit natural connections from the Levi-Civita connection of the ambient manifold.  They also play a central role in the Atiyah--Drinfeld--Hitchin--Manin (ADHM) construction of instantons, which constructs solutions of the anti-self-dual Yang-Mills equation using induced connections (for reviews, see \cite{CWS,donkrom,jardim,mantonsutcliffe}).

A fundamental property of any connection is the collection of its parallel transport operators.  These are linear maps between fibres $E_p\to E_q$ that depend on a choice of curve $\gamma$ from $p$ to $q$.  The main result of this note is a collection of high-order numerical methods to approximate the parallel transport operators of an induced connection.

There is already an obvious method to approximate parallel transport of an induced connection.  Computing parallel transport amounts to solving an initial value problem
\begin{equation}
\label{parallel transport 0}
\Omega'(x,0) = - A(x)\Omega(x,0), \quad \Omega(0,0)=\mathrm{Id},
\end{equation}
in which $A=v^\dagger\nabla_{\gamma'}v$ is the connection matrix with respect to a chosen orthonormal frame $v(x)$.  So, one may approximate $A$ using finite differences, and approximate the solution $y$ using a Runge-Kutta method, for example.  There are two problems with this naive method.  The first is that it is not gauge-covariant.  The matrix $\Omega(x,0)$ represents a linear map between the fibres of $E$ at 0 and $x$, so depends on the choice of bases $v(0)$ and $v(x)$, but does not depend on the choice of bases $v(w)$ at intermediate points $0<w<x$.  However, any approximate solution obtained using the naive method described above would depend on $v(w)$ at intermediate points $w$.  In particular, if $v(w)$ is chosen to depend on $w$ in a highly discontinuous way then the accuracy of methods such as Runge-Kutta (which assume analyticity of all functions involved) is questionable.

A second criticism of the naive method is that it is inefficient.  It entails computing derivatives (to obtain $A$) and then partially undoing this by computing integrals (to solve the parallel transport equation).  The methods that we present below compute parallel transport directly from $v(x)$, so avoid this inefficency.

The results that we presented here are motivated by ongoing work to approximate skyrmions using the Atiyah--Manton--Sutcliffe construction \cite{am1,am2,CH,CHW,HH,houghton,LM,LMS,SS,buckyball,sutcliffe,sutcliffe2}.  This construction approximates soliton solutions of a nonlinear field theory, the Skyrme model, using parallel transport operators of Yang-Mills instantons.  The simplest and most effective way of constructing instantons is the ADHM method, which uses induced connections.  To obtain a skyrmion from an instanton entails computing hundreds of parallel transports.  Moreover, to compute quantities relevant to applications in nuclear physics one may need to compute skyrmions from hundreds of different instantons.  So an efficient and accurate method to compute parallel transport of an induced connection is highly desirable in this context.

Given the ubiquity of induced connections, it seems likely that our results will prove useful in other contexts.  We sketch one possible further application to the geometry of curves at the end of section \ref{sec:4}.  An outline of this paper is as follows: in section \ref{sec:2} we establish our notation and derive some simple approximations to parallel transport.  In section \ref{sec:3} we introduce an operator formalism and use this to derive more sophisticated approximations to parallel transport.  In section \ref{sec:4} we illustrate our method in a simple example and describe some applications.  Section \ref{sec:5} discusses some interesting theoretical questions about our method.

\section{Simple approximations to parallel transport}
\label{sec:2}

\subsection{Statement of the problem}
Throughout this article we will take our ambient vector bundle to be the trivial bundle $\mathbb{C}^n\times \mathbb{R}$ over the manifold $\mathbb{R}$, equipped with the standard hermitian metric and the trivial connection.  No generality is lost here, because parallel transport is always defined along a 1-dimensional submanifold of the ambient manifold, and all vector bundles and connections over $\mathbb{R}$ are trivial.

We will let $E$ be a rank $m$ sub-bundle with orthonormal frame $v^1,\ldots v^m$.  This means that, for each $x\in\mathbb{R}$, the fibre $E_x$ is the span of vectors $v^1(x),\ldots v^m(x)$ satisfying $v^i(x)^\dagger v^j(x)=\delta_{ij}$.  Let $v$ be the $n\times m$ matrix-valued function whose columns are $v^1,\ldots v^m$; note that $v^\dagger v=\mathrm{Id}_m$.  A section of $E$ can be written in the form
\begin{equation}
z(x) = v(x)y(x),\quad y:\mathbb{R}\to \mathbb{C}^m.
\end{equation}
A section $z$ is parallel if $z'\in E^\perp$.  This is equivalent to $v^\dagger z'=0$, which is in turn equivalent to
\begin{equation}\label{parallel transport}
y'(x) + A(x)y(x)=0,\quad A(x):=v(x)^\dagger v'(x).
\end{equation}
Equation \eqref{parallel transport} is known as the parallel transport equation, and $A=v^\dagger v'$ is the matrix of the induced connection.

The solution of the parallel transport equation with initial condition $y(x_0)=y_0$ can be written
\begin{equation}
y(x)=\Omega(x,x_0)y_0,
\end{equation}
where $\Omega$ is a $U(m)$-valued function, called the parallel transport operator.  The parallel transport operator is the unique solution to
\begin{equation}
\label{Omega ODE}
\frac{d}{dx}\Omega(x,x_0)+A(x)\Omega(x,x_0)=0,\quad \Omega(x_0,x_0)=\mathrm{Id}_m.
\end{equation}

Gauge transformations $g:\mathbb{R}\to U(m)$ correspond to an $x$-dependent change of basis.  They act as
\begin{equation}
v(x) \mapsto v(x)g(x),\quad y(x)\mapsto g(x)^\dagger y(x).
\end{equation}
The induced action on the parallel transport operator is
\begin{equation}
\Omega(x,x_0) \mapsto g(x)^\dagger\Omega(x,x_0)g(x_0).
\end{equation}

The goal of this article is to find approximations $\Omega^k(x+h,x)$ to $\Omega(x+h,x)$, such that
\begin{equation}
\Omega(x+h,x)=\Omega^k(x+h,x)+O(h^{k+1}).
\end{equation}
Our approximations will be written as rational functions of $v(x_i)$ for a finite set of points $x_0<x_1<x_2<\ldots$.  We will require that, under gauge transformations $v(x_i)\mapsto v(x_i)g(x_i)$, $\Omega^k$ transforms in the same way as $\Omega$:
\begin{equation}
\Omega^k(x+h,x)\mapsto g(x+h)^\dagger\Omega^k(x+h,x) g(x).
\end{equation}

\subsection{Order 2 approximation}

A simple solution to this problem (used earlier in \cite{CHW}) is
\begin{equation}
\label{order 1 approximation}
\Omega^1(x+h,x)=v(x+h)^\dagger v(x).
\end{equation}
Notice that under gauge transformations, $v(x+h)^\dagger v(x)\mapsto g(x+h)^\dagger v(x+h)^\dagger v(x)g(x)$, so $\Omega^1$ transforms in the desired way.

To see that the approximation is order 1, we use Taylor expansions.  The parallel transport equation \eqref{parallel transport} implies that
\begin{align}
y' &= -v^\dagger v'y \\
y'' &= - (v')^\dagger v' y -v^\dagger v'' y - v^\dagger v' y'\\
&= ( -(v')^\dagger v' -v^\dagger v''+(v^\dagger v')^2 )y .
\end{align}
So
\begin{align}
\Omega(x+h,x)y(x)&=y(x+h)\\
&=y(x_0)+hy'(x)+\frac{h^2}{2}y''(x)+O(h^3)\\
%&= \left[ 1 -hv(x)^\dagger v'(x) + \frac{h^2}{2}( -(v'(x)^\dagger)v'(x) -v(x)^\dagger v''(x)+(v(x)^\dagger v'(x))^2 )\right]y(x) + O(h^3)\\
&= \left[ 1 -hv^\dagger v' + \frac{h^2}{2}( -(v')^\dagger v' -v^\dagger v''+(v^\dagger v')^2 )\right]y + O(h^3),\label{Omega order 2 expansion}
\end{align}
where in the last line $y,v,v',v''$ are understood to be evaluated at $x$.  On the other hand,
\begin{align}
v(x+h)^\dagger v(x)&= \left[v(x)+hv'(x)+\frac{h^2}{2}v''(x)+O(h^3)\right]^\dagger v(x) \\
&= v^\dagger v + h (v')^\dagger v +\frac{h^2}{2}(v'')^\dagger v+O(h^3)\\
&= 1 - h v^\dagger v' -h^2\left((v')^\dagger v'+\frac{1}{2}v^\dagger v''\right)+O(h^3),
\end{align}
where in the last line we used $v^\dagger v=1$, $(v')^\dagger v+v^\dagger v'=0$ and $(v'')^\dagger v + 2(v')^\dagger v'+v^\dagger v''=0$.  Comparing the two calculations, we see that $\Omega(x+h,x)=v^\dagger(x+h)v(x)+O(h^2)$.

In order to improve this method, we seek a second order approximation in the form
\begin{equation}
\Omega^2(x+h,x)=a\,v(x+h)^\dagger v(x)+b[v(x)^\dagger v(x+h)]^{-1}
\end{equation}
where $a,b\in\mathbb{R}$ are to be determined.  Note that the choice of operators on the right ensures that $\Omega^2$ transforms in the desired way under gauge transformations.  To compare this with $\Omega$ we need the Taylor expansion of the second operator:
\begin{align}
[v(x)^\dagger v(x+h)]^{-1} &= \left[v^\dagger v+hv^\dagger v'+\frac{h^2}{2}v^\dagger v''+O(h^3)\right]^{-1}\\
&= 1 - hv^\dagger v' + h^2\left((v^\dagger v')^2-\frac{1}{2}v^\dagger v''\right)+O(h^3).
\end{align}
So our approximation is
\begin{multline}
a\,v(x+h)^\dagger v(x)+b[v(x)^\dagger v(x+h)]^{-1} = \\
(a+b)(1-hv^\dagger v') + h^2\left(-\frac{a+b}{2}v^\dagger v''-a(v')^\dagger v'+b(v^\dagger v')^2\right)+O(h^3).
\end{multline}
This agrees with the expansion \eqref{Omega order 2 expansion} of $\Omega$ precisely when $a=b=\frac{1}{2}$.  So our order two method is
\begin{equation}
\label{Omega2}
\Omega^2(x+h,x)=\frac{1}{2}\left(v(x+h)^\dagger v(x)+[v(x)^\dagger v(x+h)]^{-1}\right).
\end{equation}

\section{Higher order approximations}
\label{sec:3}

\subsection{An operator expression for $\Omega$}
The method used above to derive a second order approximation can in principle be used to derive higher order approximations.  However, in practice the algebra quickly becomes cumbersome.  In this section we derive an operator expression for $\Omega$ that allows for much simpler derivation of approximations $\Omega^k$.

Recall that a section of $E$ is a function $z:\mathbb{R}\to\mathbb{C}^m$ such that $vv^\dagger z=z$.
Consider the operator $\omega(h)$ acting on such sections $z$ as follows:
\begin{equation}
\label{omega and Omega}
(\omega(h)z)(x) = v(x)\Omega(x,x-h)v^\dagger(x-h)z(x-h).
\end{equation}
The matrices $\Omega(x+h,x)$ determine, and are determined by, the operators $\omega(h)$.  The advantage of introducing the operators $\omega(h)$ is that they can be expressed in the following simple way:
\begin{equation}
\label{omega and A}
\omega(h)z = v\exp(-h(d+A))v^\dagger z.% = v\exp(-hv^\dagger d v)v^\dagger z.
\end{equation}
In this expression, $d$ denotes the operator $d\,z=z'$, and $v,v^\dagger,A$ act on vector-valued functions by matrix multiplication.

To show that the right hand sides of \eqref{omega and A} and \eqref{omega and Omega} are equal, we first consider the case where $h=0$.  In this situation both \eqref{omega and A} and \eqref{omega and Omega} correspond to multiplying $z(x)$ with the identity matrix, so they agree.

To show that they agree for all values of $h$, we differentiate both:
\begin{align}
\left(\frac{\partial}{\partial x}+\frac{\partial}{\partial h}\right)&v(x)\Omega(x,x-h)v^\dagger(x-h)z(x-h)\\
&= \big(v'(x)-v(x)A(x)\big)\Omega(x,x-h)v^\dagger(x-h)z(x-h)\\
&= \big(v'v^\dagger-vAv^\dagger\big) v(x)\Omega(x,x-h)v^\dagger(x-h)z(x-h)\\
\left(\frac{\partial}{\partial x}+\frac{\partial}{\partial h}\right)&(v\exp(-h(d+A))v^\dagger z)\\
 &= d(v\exp(-h(d+A))v^\dagger z)-v(d+A)\exp(-h(d+A))v^\dagger z \\
&= \big([d,v]-vA\big)\exp(-h(d+A))v^\dagger z \\
&= \big(v'v^\dagger-vAv^\dagger\big)v\exp(-h(d+A))v^\dagger z.
\end{align}
In both cases, we find that $\frac{\partial}{\partial h}\omega(h)=(-d-vAv^\dagger+v'v^\dagger)\omega(h)$.  Since the two operators satisfy the same differential equation, they agree for all values of $h$.

It will prove convenient to rewrite \eqref{omega and A} as follows.  Since $v^\dagger v$ is the identity operator and $A=v^\dagger[d,v]$, we have that $d+A=v^\dagger d\,v$.  Therefore
\begin{align}
\omega(h) &= v\exp(-hv^\dagger d\,v) v^\dagger\\
&= p-h\,p\,d\,p+\tfrac{h^2}{2}p\,d\,p\,d\,p - \tfrac{h^3}{6}p\,d\,p\,d\,p\,d\,p+\tfrac{h^4}{24}p\,d\,p\,d\,p\,d\,p\,d\,p+O(h^5),\label{omega expansion}
\end{align}
in which $p=vv^\dagger$.
We now seek operators $\omega^k(h)$ that approximate $\omega(h)$, and from these deduce approximations to $\Omega(x+h,x)$.

\subsection{Order 2}

We begin by rederiving the second order expression obtained earlier.  We seek an approximation $\omega^2(h)$ to $\omega(h)$ using the operators $\pi_1,\pi_2$ defined by
\begin{align}
(\pi_1(h)z)(x) &= v(x)v^\dagger(x-h)z(x-h)\\
\pi_2(h) &= \pi_1(-h)^{-1}.
\end{align}
We obtain Taylor expansions as follows:
\begin{align}
\pi_1(h) &= v \exp(-dh) v^\dagger\\
&= p-h\,p\,d\,p+\tfrac{h^2}{2}p\,d^2p+O(h^3)\label{pi1}\\
\pi_2(h) &= (p+h\,p\,d\,p+\tfrac{h^2}{2}p\,d^2p)^{-1}+O(h^3)\\
&= p-h\,p\,d\,p-\tfrac{h^2}{2}p\,d^2p+h^2p\,d\,p\,d\,p+O(h^3)\label{pi2}.
\end{align}
The expansion for $\pi_2$ was derived from the identities $\pi_1(-h)\pi_2(h)=p$ and $p^2=p$, together with the fact that $p \,z=z$ for any section $z$ of the subbundle $E$.

If we compare equations \eqref{pi1} and \eqref{pi2} with the expansion \eqref{omega expansion} of $\omega(h)$ it is clear that
\begin{equation}
\omega^2(h) = \frac{1}{2}(\pi_1(h)+\pi_2(h))
\end{equation}
is a second order approximation to $\omega(h)$.  The corresponding approximation to $\Omega$ is that given in \eqref{Omega2}

\subsection{Order 3}
To obtain an order 3 approximation, we consider an ansatz
\begin{equation}\label{omega3 ansatz}
\omega^3(h)=a_1\,\omega^2(\tfrac{h}{2})\omega^2(\tfrac{h}{2})+a_2\,\omega^2(h).
\end{equation}
To compare this with \eqref{omega expansion} we need an expansion for $\omega^2$.  This is obtained in a similar way to the expansions \eqref{pi1} and \eqref{pi2}:
\begin{multline}
\omega^2(h)
=p-h\,p\,d\,p+\tfrac{h^2}{2}p\,d\,p\,d\,p\\
+h^3\big(-\tfrac{1}{6}p\,d^3p+\tfrac{1}{4}p\,d^2p\,d\,p+\tfrac{1}{4}p\,d\,p\,d^2p-\tfrac{1}{2}p\,d\,p\,d\,p\,d\,p\big)\\
+h^4\big(\tfrac{1}{12}p\,d^3p\,d\,p+\tfrac{1}{8}p\,d^2p\,d^2p+\tfrac{1}{12}p\,d\,p\,d^3p\\
-\tfrac{1}{4}p\,d^2p\,d\,p-\tfrac{1}{4}p\,d\,p\,d^2p\,d\,p-\tfrac{1}{4}p\,d\,p\,d\,p\,d^2p
+\tfrac{1}{2}p\,d\,p\,d\,p\,d\,p\,d\,p\big)+O(h^5).\label{omega2 expansion}
\end{multline}
From this it follows that
\begin{multline}
\omega^2(\tfrac{h}{2})\omega^2(\tfrac{h}{2}) = p-h\,p\,d\,p+\tfrac{h^2}{2}p\,d\,p\,d\,p\\+h^3(-\tfrac{1}{24}p\,d^3p+\tfrac{1}{16}p\,d\,p\,d^2p+\tfrac{1}{16}p\,d^2p\,d\,p-\tfrac{1}{4}p\,d\,p\,d\,p\,d\,p)+O(h^4).
\end{multline}
By comparing these two expansions with \eqref{omega expansion}, we see that $\omega^3$ given in \eqref{omega3 ansatz} is an order 3 approximation to $\omega$ if and only if
\begin{equation}\label{order 3 linear system}
\begin{pmatrix}
1&1\\
-1&-1\\
\frac12&\frac12\\
-\frac{1}{24}&-\frac{1}{6}\\
\frac{1}{16}&\frac{1}{4}\\
\frac{1}{16}&\frac{1}{4}\\
-\frac{1}{4}&-\frac{1}{2}
\end{pmatrix}
\begin{pmatrix}a_1\\a_2\end{pmatrix}=
\begin{pmatrix}1\\-1\\\frac{1}{2}\\0\\0\\0\\-\frac{1}{6}\end{pmatrix}.
\end{equation}
The unique solution is $a_1=\frac43$ and $a_2=-\frac13$.  The corresponding order 3 approximation to $\omega$ is
\begin{equation}
\Omega^3(x+h,x)=\frac{4}{3}\Omega^2(x+h,x+\tfrac{h}{2})\Omega^2(x+\tfrac{h}{2},x)-\frac{1}{3}\Omega^2(x+h,x).
\end{equation}

%By comparing these two expansions with \eqref{omega expansion}, we see that $\omega^3$ given in \eqref{omega3 ansatz} is an order 3 approximation to $\omega$ if and only if $a_1+a_2=1$ and $\frac{1}{4}a_1+\frac{1}{2}a_2=\frac{1}{6}$.  The unique solution is $a_1=\frac43$ and $a_2=-\frac13$.  The corresponding order 3 approximation to $\omega$ is
%\begin{equation}
%\Omega^3(x+h,x)=\frac{4}{3}\Omega^2(x+h,x+\tfrac{h}{2})\Omega^2(x+\tfrac{h}{2},x)-\frac{1}{3}\Omega^2(x+h,x).
%\end{equation}

\subsection{Order 4}
To obtain an order 4 approximation, we consider operators of the form
\begin{equation}\label{omega4 ansatz}
\omega^4(h) = a_1\omega^2(\tfrac{h}{3})\omega^2(\tfrac{h}{3})\omega^2(\tfrac{h}{3})+a_2\omega^2(\tfrac{h}{3})\omega^2(\tfrac{2h}{3})+a_3\omega^2(\tfrac{2h}{3})\omega^2(\tfrac{h}{3})+a_4\omega^2(h).
\end{equation}
The expansion of the final operator appearing on the right is given in \eqref{omega2 expansion}, and the expansions of the remaining three operators can all be derived from \eqref{omega2 expansion}:
\begin{multline}
\omega^2(\tfrac{h}{3})\omega^2(\tfrac{h}{3})\omega^2(\tfrac{h}{3})
= p-h\,p\,d\,p+\tfrac{h^2}{2}p\,d\,p\,d\,p\\
+h^3\big(-\tfrac{1}{54}p\,d^3p+\tfrac{1}{36}p\,d^2p\,d\,p+\tfrac{1}{36}p\,d\,p\,d^2p-\tfrac{11}{54}p\,d\,p\,d\,p\,d\,p\big)\\
+h^4\big(\tfrac{1}{108}p\,d^3p\,d\,p+\tfrac{1}{216}p\,d^2p\,d^2p+\tfrac{1}{108}p\,d\,p\,d^3p\\
-\tfrac{1}{54}p\,d^2p\,d\,p\,d\,p-\tfrac{1}{36}p\,d\,p\,d^2p\,d\,p-\tfrac{1}{54}p\,d\,p\,d\,p\,d^2p
+\tfrac{1}{12}p\,d\,p\,d\,p\,d\,p\,d\,p\big)+O(h^5),
\end{multline}
\begin{multline}
%\frac{1}{4}v\left(v^\dagger e^{-\tfrac{h}{3}d}v+\left[v^\dagger e^{\tfrac{h}{3}d}v\right]^{-1}\right)\left(v^\dagger e^{-\tfrac{2h}{3}d}v+\left[v^\dagger e^{\tfrac{2h}{3}d}v\right]^{-1}\right)v^\dagger
\omega^2(\tfrac{h}{3})\omega^2(\tfrac{2h}{3})
= p-h\,p\,d\,p+\tfrac{h^2}{2}p\,d\,p\,d\,p\\
+h^3\big(-\tfrac{1}{18}p\,d^3p+\tfrac{1}{12}p\,d^2p\,d\,p+\tfrac{1}{12}p\,d\,p\,d^2p-\tfrac{5}{18}p\,d\,p\,d\,p\,d\,p\big)\\
+h^4\big(\tfrac{7}{324}p\,d^3p\,d\,p+\tfrac{17}{648}p\,d^2p\,d^2p+\tfrac{11}{324}p\,d\,p\,d^3p\\
-\tfrac{19}{324}p\,d^2p\,d\,p\,d\,p-\tfrac{1}{12}p\,d\,p\,d^2p\,d\,p-\tfrac{25}{324}p\,d\,p\,d\,p\,d^2p
+\tfrac{29}{162}p\,d\,p\,d\,p\,d\,p\,d\,p\big)+O(h^5),
\end{multline}
\begin{multline}
%\frac{1}{4}v\left(v^\dagger e^{-\tfrac{2h}{3}d}v+\left[v^\dagger e^{\tfrac{2h}{3}d}v\right]^{-1}\right)\left(v^\dagger e^{-\tfrac{h}{3}d}v+\left[v^\dagger e^{\tfrac{h}{3}d}v\right]^{-1}\right)v^\dagger\\
\omega^2(\tfrac{2h}{3})\omega^2(\tfrac{h}{3})
= p-h\,p\,d\,p+\tfrac{h^2}{2}p\,d\,p\,d\,p\\
+h^3\big(-\tfrac{1}{18}p\,d^3p+\tfrac{1}{12}p\,d^2p\,d\,p+\tfrac{1}{12}p\,d\,p\,d^2p-\tfrac{5}{18}p\,d\,p\,d\,p\,d\,p\big)\\
+h^4\big(\tfrac{11}{324}p\,d^3p\,d\,p+\tfrac{17}{648}p\,d^2p\,d^2p+\tfrac{7}{324}p\,d\,p\,d^3p\\
-\tfrac{25}{324}p\,d^2p\,d\,p\,d\,p-\tfrac{1}{12}p\,d\,p\,d^2p\,d\,p-\tfrac{19}{324}p\,d\,p\,d\,p\,d^2p
+\tfrac{29}{162}p\,d\,p\,d\,p\,d\,p\,d\,p\big)+O(h^5).
\end{multline}
It follows that $\omega^4$ given in equation \eqref{omega4 ansatz} is an order 4 approximation to $\omega$ given in \eqref{omega expansion} if and only if $a_1,a_2,a_3,a_4$ satisfy the linear equation,
\begin{equation}\label{order 4 linear system}
\begin{pmatrix}
1&1&1&1\\
-1&-1&-1&-1\\
%0&0&0&0\\
\frac{1}{2}&\frac{1}{2}&\frac{1}{2}&\frac{1}{2}\\
-\frac{1}{54}&-\frac{1}{18}&-\frac{1}{18}&-\frac{1}{6}\\
\frac{1}{36}&\frac{1}{12}&\frac{1}{12}&\frac{1}{4}\\
\frac{1}{36}&\frac{1}{12}&\frac{1}{12}&\frac{1}{4}\\
-\frac{11}{54}&-\frac{5}{18}&-\frac{5}{18}&-\frac{1}{2}\\
%0&0&0&0\\
\frac{1}{108}&\frac{7}{324}&\frac{11}{324}&\frac{1}{12}\\
\frac{1}{216}&\frac{17}{648}&\frac{17}{648}&\frac{1}{8}\\
\frac{1}{108}&\frac{11}{324}&\frac{7}{324}&\frac{1}{12}\\
-\frac{1}{54}&-\frac{19}{324}&-\frac{25}{324}&-\frac{1}{4}\\
-\frac{1}{36}&-\frac{1}{12}&-\frac{1}{12}&-\frac{1}{4}\\
-\frac{1}{54}&-\frac{25}{324}&-\frac{19}{324}&-\frac{1}{4}\\
\frac{1}{12}&\frac{29}{162}&\frac{29}{162}&\frac{1}{2}
\end{pmatrix}
\begin{pmatrix}a_1\\a_2\\a_3\\a_4\end{pmatrix}
=\begin{pmatrix}1\\-1\\\frac{1}{2}\\0\\0\\0\\-\frac{1}{6}\\0\\0\\0\\0\\0\\0\\\frac{1}{24}\end{pmatrix}
\end{equation}
The unique solution is given by $(a_1,a_2,a_3,a_4)=\frac{1}{44}(90,-27,-27,8)$.  So our order 4 approximation to $\Omega$ is
\begin{multline}\label{Omega 4 ansatz}
\Omega^4(x+h,x)=\frac{1}{44}\Big(90\Omega^2(x+h,x+\tfrac{2h}{3})\Omega^2(x+\tfrac{2h}{3},x+\tfrac{h}{3})\Omega^2(x+\tfrac{h}{3},x)\\
-27\Omega^2(x+h,x+\tfrac{2h}{3})\Omega^2(x+\tfrac{2h}{3},x)-27 \Omega^2(x+h,x+\tfrac{h}{3})\Omega^2(x+\tfrac{h}{3},x)\\
+8\Omega^2(x+h,x)\Big).
\end{multline}

\subsection{Improved methods}

The connection matrix $A=v^\dagger v'$ that appears in the parallel transport equation \eqref{parallel transport} is anti-hermitian.  It follows that $\Omega(x+h,x)$ is a unitary $m\times m$ matrix, and hence that
\begin{equation}
|\det(\Omega(x+h,x))| = 1.
\end{equation}
Therefore
\begin{equation}
\widehat{\Omega}^k(x+h,x):=|\det(\Omega^k(x+h,x))|^{-\frac{1}{m}}\Omega^k(x+h,x)
\end{equation}
satisfies $|\det(\widehat{\Omega}^k(x+h,x))|=1$ and hence is a better approximation to parallel transport that $\Omega^k$.  In fact, in certain situations $\widehat{\Omega}^k(x+h,x)-\Omega(x+h,x)=O(h^{k+2})$, so this improved approximation is an order of magnitude better than $\Omega^k$.  We review these situations below.

The first case to consider is where our sub-bundle $E$ has rank 1.  In this case $\Omega^k$ is a $1\times 1$ matrix, and so
\begin{equation}\label{rank1det}
\det(\Omega^k(x+h,x))=\Omega^k(x+h,x).
\end{equation}
Since $\Omega(x+h,x)$ is unitary and $\Omega(x+h,x)^{-1}=\Omega(x,x+h)$, we have that
\begin{equation}
\Omega(x+h,x)^\dagger=\Omega(x,x+h)
\end{equation}
We will assume similarly that
\begin{equation}\label{OmegaN hermitian}
\Omega^k(x+h,x)^\dagger=\Omega^k(x,x+h).
\end{equation}
This assumption is satisfied by all of the approximations $\Omega^1$, $\Omega^2$, $\Omega^3$, $\Omega^4$ obtained above.
%It follows that
%\begin{equation}
%\big(e^{-hd}\Omega^k(x+h,x)\big)^\dagger=\Omega^k(x,x+h)e^{hd}=e^{hd}\Omega^k(x-h,x).
%\end{equation}
Now we introduce $\epsilon^k(x)$ such that
\begin{equation}
\Omega^k(x+h,x)=\Omega(x+h,x)+h^{k+1}\epsilon^k(x)+O(h^{k+2}).
\end{equation}
It follows that
\begin{align}
\Omega^k(x,x+h) &= \Omega(x,x+h)+(-h)^{k+1}\epsilon^k(x+h)+O(h^{k+2}) \\
&= \Omega(x,x+h)+(-h)^{k+1}\epsilon^k(x)+O(h^{k+2}),
\end{align}
because $\epsilon^k(x+h)=\epsilon^k(x)+O(h)$.
%It follows immediately that $\Omega^k(x,x+h)=\Omega(x,x+h)-h^{k+1}\epsilon^k+O(h^{k+2})$, so
%\begin{equation}
%\epsilon^k(x+h,x)^\dagger = (-h)^{k+1}\epsilon^k(x,x+h).
%\end{equation}
Therefore
\begin{align}
|\det(\Omega^k(x+h,x))|^2 &=
\Omega^k(x+h,x)\Omega^k(x+h,x)^\dagger \hspace{30pt}\text{by \eqref{rank1det}} \label{eq:detsq}\\
&= \Omega^k(x+h,x)\Omega^k(x,x+h)\hspace{34pt}\text{by \eqref{OmegaN hermitian}}\\
&= \big(\Omega(x+h,x)+h^{k+1}\epsilon^k(x)\big)\nonumber\\
&\quad\times \big(\Omega(x,x+h)+(-h)^{k+1}\epsilon^k(x)\big)+O(h^{k+2})\\
&= 1 + \big[1+(-1)^{k+1}\big]h^{k+1}\epsilon^k + O(h^{k+2}).
\end{align}
Here in the final line we used that $\Omega(x+h,x)\Omega(x,x+h)=1$ and that $\Omega(x+h,x)=1+O(h)$.

Thus, in the case that $k$ is odd, we obtain
\begin{align}
\widehat{\Omega}^k(x+h,x)&=\Omega^k(x+h,x)\big(1 + 2h^{k+1}\epsilon^k\big)^{-\frac{1}{2}}+O(h^{k+2}) \\
&=\big(\Omega(x+h,x)+h^{k+1}\epsilon^k\big) \big(1 -h^{k+1}\epsilon^k\big) + O(h^{k+2})\\
&=\Omega(x+h,x)+O(h^{k+2})
\end{align}
Thus, if $k$ is odd and $E$ is a complex line bundle, $\widehat{\Omega}^{k}$ is an order $k+1$ approximation to $\Omega$.

We obtain a similar result if $E$ is a symplectic bundle of rank 2.  Recall that, if $n$ is even, $\mathbb{C}^n$ carries a symplectic structure defined by an anti-linear map $J:\mathbb{C}^n\to\mathbb{C}^n$ of the form
\begin{equation}
J:\begin{pmatrix}u_1\\u_2\\\vdots\\u_{n-1}\\u_{n}\end{pmatrix}\mapsto \begin{pmatrix}\bar{u}_2\\-\bar{u}_1\\\vdots\\ \bar{u}_n \\ -\bar{u}_{n-1}\end{pmatrix}.
\end{equation}
A subspace $E$ is called symplectic if $Ju\in E$ for all $u\in E$.  If $E$ is symplectic and of rank 2, it admits an orthonormal basis of the form $v_1,Jv_1$.  In this case we can write the $n\times 2$ basis matrix $v=(v_1\,Jv_1)$ in the form
\begin{equation}
v=\begin{pmatrix}q_1^0\mathbf{1}+q_1^1\mathbf{i}+q_1^2\mathbf{j}+q_1^3\mathbf{k}\\\vdots\\q_{n/2}^0\mathbf{1}+q_{n/2}^1\mathbf{i}+q_{n/2}^2\mathbf{j}+q_{n/2}^3\mathbf{k}\end{pmatrix},
\end{equation}
in which $q^\mu_i$ are real and
\begin{equation}
\mathbf{1}=\begin{pmatrix}1&0\\0&1\end{pmatrix},\,
\mathbf{i}=\begin{pmatrix}0&-\imath\\-\imath&0\end{pmatrix},\,
\mathbf{j}=\begin{pmatrix}0&-1\\1&0\end{pmatrix},\,
\mathbf{k}=\begin{pmatrix}-\imath&0\\0&\imath\end{pmatrix}.
\end{equation}
In other words, $v$ can be written as a vector of quaternions.  It follows that each of the approximations $\Omega^k(x+h,x)$ can be written as a real linear combination of $\mathbf{1},\mathbf{i},\mathbf{j},\mathbf{k}$, and hence that
\begin{equation}
|\det\Omega^k(x+h,x)|^2\mathbf{1}=\Omega^k(x+h,x)\Omega^{k}(x+h,x)^\dagger.
\end{equation}
This means that the calculation starting with equation \eqref{eq:detsq} goes through as in the rank 1 case, and we again obtain that if $k$ is odd then $\widehat{\Omega}^{k}$ is an order $k+1$ approximation to $\Omega$.

The final situation to consider is where $E$ is a real rank 1 sub-bundle of a real vector bundle.  This case is trivial in the sense that the parallel transport operator is a $1\times 1$ orthogonal matrix, so is either $1$ or $-1$.  Similarly, $\widehat{\Omega}^k(x+h,x)$ has determinant $\pm1$ so is either $1$ or $-1$.  Thus, for sufficiently small $h$, $\widehat{\Omega}^k$ is a perfect approximation to $\Omega^k$.

\section{Implementation}
\label{sec:4}
\subsection{Simple example}
We now illustrate the methods developed above in a simple example.  For $t$ in the interval $[-\pi/2,\pi/2]$, let $E_t\subset \mathbb{C}^4$ be the kernel of the matrix,
\begin{equation}
\label{ADHM matrix}
\Gamma(t)=\begin{pmatrix}\cos t &0&\sin t&\cos t\\0&\cos t&-\cos t & \sin t\end{pmatrix}.
\end{equation}
Then $E$ is a rank 2 sub-bundle of the trivial rank 4 bundle over $[-\pi/2,\pi/2]$.  We will approximate the parallel transport operator $\Omega(\pi/2,-\pi/2)$.

To do so, we choose $N+1$ equally-spaced points $t_i=-\pi/2+i\pi/N$ in the interval $[-\pi/2,\pi/2]$, with $0\leq i\leq N$.  For each $i$ we find an orthonormal basis for the kernel of $\Gamma(t)$ and arrange the basis vectors into a $4\times 2$ matrix $v(t_i)$ satisfying $v(t_i)^\dagger v(t_i)=\mathrm{Id}_2$.  The kernels of $\Gamma(t_0)$ and $\Gamma(t_N)$ are equal, and for both of these we choose the basis
\begin{equation}
v(t_0)=v(t_N)=\begin{pmatrix}1&0\\0&1\\0&0\\0&0\end{pmatrix}
\end{equation}
To approximate parallel transport to accuracy $1/N^k$, we compute the matrices $\Omega^k(t_{(i+1)(k-1)},t_{i(k-1)})$ for $0\leq i<N/(k-1)$.  We then compute
\begin{equation}
\label{numerical approximation}
U = \Omega^k(t_N,t_{N-k+1)})\Omega^k(t_{N-k+1},t_{N-2k+2})\ldots \Omega^k(t_{k-1},t_0).
\end{equation}
Our earlier results imply that $\Omega(\pi/2,-\pi/2)=U+O(1/N^k)$.
%\begin{equation}
%\label{numerical approximation}
%\Omega(\pi/2,-\pi/2)= \Omega^m(t_N,t_{N-m})\Omega^m(t_{N-m},t_{N-2m})\ldots \Omega^m(t_m,t_0) + O\left(\frac{1}{N^m}\right).
%\end{equation}

The matrix $\Omega(\pi/2,-\pi/2)$ can in fact be computed exactly by solving the differential equation \eqref{Omega ODE}.  The result is
\begin{equation}
\Omega(\pi/2,-\pi/2)=\begin{pmatrix}-\cos(\pi/\sqrt{2})&-\sin(\pi/\sqrt{2})\\\sin(\pi/\sqrt{2})&-\cos(\pi/\sqrt{2})\end{pmatrix} .
\end{equation}
We can assess the accuracy of our approximation by computing the error $E=\frac{1}{2}\operatorname{Tr}(\Delta\Delta^\dagger)$, in which $\Delta=U-\Omega(\pi/2,-\pi/2)$.  The results are displayed in table \ref{table1}.  As expected, using a higher order method allows one to attain a desired accuracy with fewer points than would be necessary with a lower order method.
\begin{table}
\begin{tabular}{r|cccccc}
$N$ & 6 & 12 & 24 & 48 & 96\\% & 192 \\
\hline
$\Omega^1$ & $3.56\times10^{-1}$ & $1.96\times10^{-1}$ & $1.03\times10^{-1}$ & $5.31\times10^{-2}$ & $2.69\times10^{-2}$ &\\% $1.35\times10^{-2}$\\
$\Omega^2$ & $1.25\times 10^{-1}$ & $3.04\times 10^{-2}$ & $7.55\times 10^{-3}$ & $1.88\times 10^{-3}$ & $4.71\times10^{-4}$\\% & $1.18\times 10^{-4}$ \\
$\Omega^3$ & $5.61\times 10^{-3}$ & $2.28\times 10^{-3}$ & $3.78\times 10^{-4}$ & $5.01\times10^{-5}$ & $6.36\times 10^{-6}$ &\\% $7.97\times10^{-7}$\\
$\Omega^4$ & $2.79\times 10^{-2}$ & $2.23\times10^{-3}$ & $7.94\times 10^{-5}$ & $3.08\times10^{-6}$ & $1.51\times10^{-7}$ &\\% $8.70\times10^{-9}$
\end{tabular}
\caption{Error in calculating the parallel transport using the methods $\Omega^k$ over $N$ intervals.}
\label{table1}
\end{table}

Now we consider the improved method.  Recall that the improved method asks as to multiply each matrix $\Omega^k(t_{(i+1)(k-1)},t_{i(k-1)})$ with a positive real number so that the modulus of its determinant is 1.  Since scalar multiplication commutes with matrix multiplication and determinants are multiplicative, this is equivalent to computing $U$ as in \eqref{numerical approximation} and then computing $\widehat{U}=U/\sqrt{\det(U)}$.  Thus the additional computational cost associated with the improved method is minimal.

Nevertheless, in the cases where $k$ is odd the improved method $\widehat{\Omega}^k$ is a substantial improvement over $\Omega^k$ and comparable in accuracy with $\widehat{\Omega}^{k+1}$, as can be seen in the table \ref{table2}.  The reason for this improvement is that the kernel of our matrix $E$ is a rank 2 symplectic subspace of $\mathbb{C}^4$, so by our earlier results $\widehat{U}-\Omega(\pi/2,-\pi/2)=O(1/N^{k+1})$ when $k$ is odd (whereas $U-\Omega(\pi/2,-\pi/2)=O(1/N^{k})$).

Table \ref{table2} also shows that the errors obtained with the methods $\widehat{\Omega}^1$ and $\widehat{\Omega}^2$ are identical: this is because these two methods are in fact mathematically equivalent.  To see this, one simply needs to note that
\begin{equation}
(v^{\dagger}(x)v(x+h))^{-1}=\frac{(v^{\dagger}(x)v(x+h))^\dagger}{|v^{\dagger}(x)v(x+h)|^2}=\frac{v^{\dagger}(x+h)v(x)}{|v^{\dagger}(x)v(x+h)|^2}
\end{equation}
using the fact that $v^\dagger(x)v(x+h)$ can be written as a quaternion.  Thus $\Omega^1$ and $\Omega^2$ given in equations \eqref{order 1 approximation} and \eqref{Omega2} agree up to scalar multiplication, and their normalised counterparts $\widehat{\Omega}^1$ and $\widehat{\Omega}^2$ agree exactly.

\begin{table}
\begin{tabular}{r|cccccc}
$N$ & 6 & 12 & 24 & 48 & 96\\% & 192 \\
\hline
$\widehat{\Omega}^1$ & $1.23\times10^{-1}$ & $3.03\times10^{-3}$ & $7.54\times10^{-3}$ & $1.88\times10^{-3}$ & $4.71\times10^{-4}$ \\% & $1.18\times10^{-4}$ \\
$\widehat{\Omega}^2$ & $1.23\times10^{-1}$ & $3.03\times10^{-3}$ & $7.54\times10^{-3}$ & $1.88\times10^{-3}$ & $4.71\times10^{-4}$\\% & $1.18\times10^{-4}$ \\
$\widehat{\Omega}^3$ & $4.30\times10^{-5}$ & $6.50\times10^{-4}$ & $4.27\times10^{-5}$ & $2.64\times10^{-6}$ & $1.65\times10^{-7}$\\% & $1.03\times10^{-8}$ \\
$\widehat{\Omega}^4$ & $9.80\times10^{-3}$ & $4.39\times10^{-4}$ & $3.52\times10^{-5}$ & $2.19\times10^{-6}$ & $1.36\times10^{-7}$\\% & $8.44\times10^{-9}$
\end{tabular}
\caption{Error in calculating the parallel transport using the improved methods $\widehat{\Omega}^k$ over $N$ intervals.}
\label{table2}
\end{table}

\subsection{Application to instantons and skyrmions}

The ADHM construction produces instantons (i.e.\ finite-action solutions of the self-dual Yang-Mills equations) using induced connections \cite{ADHM,CWS,donkrom,jardim,mantonsutcliffe}.  In fact, all instantons can be produced by the ADHM method.

In the case of gauge group $SU(2)$, the method starts with an $(n+1)\times n$ matrix $\Delta$ of quaternions that depends on a point $x\in\mathbb{R}^4$.  This must be written in the form
\begin{equation}
\Delta(x^1,x^2,x^3,x^4) = \begin{pmatrix}L \\ M - (x_4\mathbf{1}+x_1\mathbf{i}+x_2\mathbf{j}+x_3\mathbf{k})\otimes \mathrm{Id}_n \end{pmatrix}
\end{equation}
with $L$, $M$ matrices of quaternions of size $1\times n$ and $n\times n$ and $M$ is symmetric.  The matrix must be such that $\Delta(x)^\dagger\Delta(x)$ is a real invertible matrix for all $x$.  This constraint ensures that the kernel $E_x$ of $\Delta(x)$ is of quaternion dimension 1 (or complex dimension 2).  So one can choose a column vector $v(x)$ of quaternions satisfying $v(x)^\dagger v(x)=\mathbf{1}$ that spans the kernel.  The instanton is obtained by setting
\begin{equation}
\label{ADHM connection}
A_\mu(x) = v(x)^\dagger \frac{\partial v}{\partial x^\mu}(x).
\end{equation}
In other words, the instanton is the induced connection on the sub-bundle $E$.

Atiyah--Manton proposed \cite{am1,am2} that holonomy of instantons could be used to approximate skyrmions, which are used to model atomic nuclei.  To be more precise, let $A$ be a fixed instanton with gauge group $SU(2)$.  For each $(x^1,x^2,x^3)\in\mathbb{R}^3$, let $U(x^1,x^2,x^3)$ be the parallel transport operator from $t=-\infty$ to $t=\infty$ along the line in $\mathbb{R}^4$ parametrised as $t\mapsto(x^1,x^2,x^3,t)$.  Atiyah--Manton proposed that the resulting function $U:\mathbb{R}^3\to SU(2)$ can be used to approximate a solution of the Euler-Lagrange equations of the Skyrme model.

This approximation was shown to work well in a number of situations \cite{houghton,LM,LMS,mantonsutcliffe,SS,buckyball}.  Subsequently, Sutcliffe gave a theoretical explanation of the success of the approximation \cite{sutcliffe}.  Sutcliffe moreover showed that instantons could also be used to approximate skyrmions coupled to vector mesons.  In Sutcliffe's construction one chooses a gauge transformation $g:\mathbb{R}^4\to SU(2)$ such that the gauge transformed connection,
\begin{equation}
\label{gauge transformation}
\tilde{A}_\mu = g^{-1}A_\mu g + g^{-1}\frac{\partial g}{\partial x^\mu},
\end{equation}
satisfies $A_4=0$.  The vector mesons are then obtained by computing the integrals,
\begin{equation}
\label{vector meson}
W_i(x^1,x^2,x^3)=\int_{-\infty}^\infty \phi(t) \tilde{A}_i(x^1,x^2,x^3,t)\,dt
\end{equation}
for a certain function $\phi(t)$.

Our methods provide an efficient numerical implementation of the Atiyah--Manton--Sutcliffe construction.  If the ADHM data of an instanton is known, then to compute the holonomy matrix $U$ at a point $(x^1,x^2,x^3)$ one needs to divide the corresponding line in $\mathbb{R}^4$ into a finite number of sub-intervals.  The holonomy matrix $U$ can then be computed as a product of parallel transport operators along these intervals.  In order to obtain accurate results it is important to include the points $(x^1,x^2,x^3,\pm\infty)$ at the two ends of the line.  The basis matrices $v(x^1,x^2,x^3,\pm\infty)$ at these two points by definition span the kernel of
\begin{equation}
\lim_{x^4\to\pm\infty}\frac{1}{x^4}\Delta(x^1,x^2,x^3,x^4)^\dagger = \pm\begin{pmatrix}0 & \mathbf{1}\otimes \mathrm{Id}_n\end{pmatrix}.
\end{equation}
We note that the kernel of this matrix is the same in both the $+\infty$ and $-\infty$ cases, and does not depend on $x^1,x^2,x^3$.  It is important to choose the same basis $v(x^1,x^2,x^3,\pm\infty)=v_\infty$ for all values of $x^1,x^2,x^3$.  A very natural choice is
\begin{equation}
v(x^1,x^2,x^3,\pm\infty) = \begin{pmatrix}\mathbf{1} \\ 0 \\ \vdots \\ 0\end{pmatrix}.
\end{equation}
The choice of bases at other points is of no consequence, because the approximation $\Omega^k$ depends only on the choice of bases at $x^4=\pm\infty$.  The only constraint is that the columns of $v$ should be orthonormal vectors in $\mathbb{C}^{2n+2}$.

The example described in the previous subsection corresponds to taking the holonomy of a charge 1 instanton.  In this case $n=1$ and the ADHM matrix is particularly simple and is given by $L=\mathbf{1}$, $M=0$.  The matrix in equation \eqref{ADHM matrix} is just $\cos t\,\Delta^\dagger(0,1,0,\tan t)$, and the parallel transport operator $\Omega(\pi/2,-\pi/2)$ that we computed was therefore $U(0,1,0)$.  Notice that our choice of parametrisation $x^4=\tan t$ maps the points $x^4=\pm\infty$ to $t=\pm\pi/2$.  This parametrisation also ensures that the points on the circle in $S^4$ that corresponds under stereographic projection to our line in $\mathbb{R}^4$ are fairly evenly spaced.  This is a sensible way to choose points, because the instanton on $\mathbb{R}^4$ constructed by the ADHM construction is the pull-back of an instanton on $S^4$.

Our method also aids the calculation of the vector mesons.  This is because the constraint $\tilde{A}_4=0$ imposed on the connection \eqref{gauge transformation} is equivalent to the parallel transport equation:
\begin{equation}
\frac{\partial g}{\partial x^4} + A_4 g = 0.
\end{equation}
One can therefore calculate $g(x^1,x^2,x^3,x^4)$ by calculating the parallel transport of $A$ along the straight line from $(x^1,x^2,x^3,-\infty)$ to $(x^1,x^2,x^3,x^4)$.  In fact, if one was also computing $U$ then one would already have calculated this parallel transport as part of that process.  Having calculated $g$, one can calculate $\tilde{A}_i$ efficiently using the identity
\begin{equation}
\tilde{A}_i = \tilde{v}^\dagger \frac{\partial \tilde{v}}{\partial x^i},\quad \tilde{v}=vg,
\end{equation}
which is easily shown to be equivalent to \eqref{ADHM connection} and \eqref{gauge transformation}.  In practice, these derivatives would be approximated as finite differences.  Finally, the integral in the definition \eqref{vector meson} of $W_i$ can be approximated by a finite sum.

\subsection{Calculating the total torsion of a space curve}

In this section we describe another possible application to the geometry of spacial curves.  Let $\mathbf{x}:[0,L]\to\mathbb{R}^3$ be a smooth arclength-parametrised closed curve (meaning that $\mathbf{x}(L)=\mathbf{x}(0)$ and $d^n\mathbf{x}/ds^{n}(L)=d^n\mathbf{x}/ds^{n}(0)$ for all $n$).  Let $\mathbf{u},\mathbf{n},\mathbf{b}$ be its Frenet frame and let $\kappa,\tau$ be its curvature and torsion.  The total torsion of $\mathbf{x}$ is
\begin{equation}
T:=\int_0^L\tau(s)\,ds.
\end{equation}
This quantity appears in a number of contexts.  For example, all curves embedded in a sphere have total torsion zero, and the sphere and plane are the only surfaces with this property \cite{millmanparker}.  The total torsion is a conserved quantity for the localised induction equation for vortex filaments \cite{hasimoto} (and is the second such quantity in the hierarchy developed in \cite{langerperline}).  The total torsion is related to the self-linking number $L\in\mathbb{Z}$ and the writhe $Wr\in\mathbb{R}$ by the formula $T/2\pi = L-Wr$ \cite{calugareanu}.

The torsion $\tau$ can be understood as the induced connection on the normal bundle to the curve.  To see this, choose the frame $v^1(s)=\mathbf{n}(s)$ and $v^2(s)=\mathbf{b}(s)$ for the normal bundle and combine these into a $3\times 2$ matrix $v$.  Then, by the Frenet equations, the induced connection is
\begin{equation}
A = v^\dagger v' = \begin{pmatrix}\mathbf{n}\cdot\mathbf{n}' & \mathbf{n}\cdot\mathbf{b}' \\ \mathbf{b}\cdot\mathbf{n}' & \mathbf{b}\cdot\mathbf{b}'\end{pmatrix} = \begin{pmatrix}0&-\tau\\\tau&0\end{pmatrix}.
\end{equation}
It follows that
\begin{equation}
\Omega(L,0) = \begin{pmatrix} \cos T & \sin T \\ -\sin T & \cos T \end{pmatrix}.
\end{equation}
Thus $\Omega(L,0)$ determines the fractional part of $T/2\pi$, and of $Wr$.

Our methods can be used to calculate the fractional part of $T/2\pi$, and hence of $Wr$, to high precision.  To do so with the order 3 method, one must first choose a finite set $\mathbf{x}_0,\ldots,\mathbf{x}_{2N-1}$ of points along the curve, written in the form $\mathbf{x}_i=\mathbf{x}(s_i)$ with $s_0<s_1<\ldots<s_{2N-1}$.  Then the curve can be approximated by a polygonal arc with edge vectors $\mathbf{u}_i:=\mathbf{x}_{i+1}-\mathbf{x}_{i}$, where the indices are understood modulo $2N$.  For each edge vector one must choose $3\times 2$ matrices $v_i$ satisfying $v_i^T \mathbf{u}_i=0$ and $v_i^Tv_i=\mathrm{Id}_2$.  One then calculates
\begin{align}
\Omega^{2}_{i}&=\frac{1}{2}(v_{i+1}^Tv_i+(v_i^Tv_{i+1})^{-1})\\
\Omega^3_{j}&=\frac{4}{3}\Omega^2_{2j+1}\Omega^2_{2j} - \frac{1}{6}(v_{2j+2}^Tv_{2j}+(v_{2j}^Tv_{2j+1})^{-1})\\
U&=\Omega^3_{N-1}\Omega^3_{N-2}\ldots\Omega^3_{0}.
\end{align}
Finally, one finds an angle $\theta\in[0,2\pi)$ by solving the system
\begin{equation}
\frac{U_{11}}{\sqrt{U_{11}^2+U_{12}^2}}=\cos\theta,\quad \frac{U_{12}}{\sqrt{U_{11}^2+U_{12}^2}}=\sin\theta.
\end{equation}
The fractional part of $T/2\pi$ is given by
\begin{equation}
\left\lfloor \frac{T}{2\pi} \right\rfloor = \frac{\theta}{2\pi}+O\left(\frac{1}{N^3}\right).
\end{equation}
Obviously, higher precision could be obtained using the order 4 method.  Note that it is not necessary to choose the basis matrices $v_i$ to approximate the Frenet frame -- any choice of orthonormal frame would be suitable, because our method respects changes of basis (and because the group $SO(2)$ is abelian).

With a little more effort, one could also compute the integer part of $T$.  To do this, one should choose the columns of $v_i$ to be discrete approximations to the normal and binormal.  For example, applying Gram--Schmidt orthogonalisation to the vectors $\mathbf{u}_i, \mathbf{u}_{i+1}-\mathbf{u}_{i-1},\mathbf{u}_i\times(\mathbf{u}_{i+1}-\mathbf{u}_{i-1})$ results in a discrete approximation to the unit tangent, normal, and binormal.  One then defines $H_j=\Omega^3_j\Omega^3_{j-1}\ldots\Omega^3_0$.  One can compute the integer part of $T/2\pi$ by looking at sign changes of the upper right entry of $H_j$.  More precisely, let $n_+$ be the number of integers $j$ such that $(H_j)_{11}>0$, $(H_j)_{12}<0$ and $(H_{j+1})_{12}\geq0$, and let $n_-$ be the number of integers $j$ such that $(H_j)_{11}>0$, $(H_j)_{12}\geq0$ and $(H_{j+1})_{12}< 0$.  Then the total torsion is given by $T=2\pi(n_+-n_-) + \theta + O(N^{-3})$.
\section{Conclusion}
\label{sec:5}
We have derived numerical methods to approximate parallel transport operators for the induced connection on a subbundle of a vector bundle.  Our methods are simpler than a naive application of the Runge--Kutta method and insensitive to choices of basis.

Our most accurate method has errors of order 4.  This level of accuracy should be sufficient for most applications.  But the algebraic framework that we have presented could be used to derive higher order methods if desired.  We expect that an order $k$ method could be obtained for any $k\in\mathbb{N}$, and it would be of interest to find a mathematical proof of (or counterexample to) this statement.

Another interesting question concerns the number of sub-intervals required to obtain an order $k$ method.  Our derivation order 3 method for approximating $\Omega(h,0)$ required us to divide the interval $[0,h]$ sub-intervals of length $h/2$, while our order 4 method requires 3 sub-intervals of length $h/3$.  We find it surprising that so few sub-intervals are needed.  To see why, one only needs to look at equation \eqref{order 4 linear system}, whose solution was required to find an order 4 method.  This is a linear system of 14 equations in 4 unknowns, so it is surprising that we were able to a solutions.  Similarly, the linear system \eqref{order 3 linear system} has a solution, despite having more equations than unknowns.  These observations suggest that there is some underlying reason why solutions can be found, despite the equations apparently being overdetermined, but we have been unable to find a satisfactory explanation.  It would be an interesting mathematical problem to determined the minimum number of subdivisions required to obtain an order $k$ method.

\bigskip
\noindent\textbf{Acknowledgement.} I am grateful to Chris Halcrow for providing feedback on the performance of these algorithms, for suggesting the improved method, and for reading a draft of this manuscript.

\bigskip
\noindent\textbf{Funding.} No funding was received to assist with the preparation of this manuscript.

\end{document}